\documentclass[12pt,oneside]{article}

\usepackage{amsfonts}
\usepackage{amscd}
\usepackage{amsthm}
\usepackage{amssymb}
\usepackage{amsmath}
\usepackage{epsfig}

\usepackage{graphicx}
\usepackage{xypic}
\input xypic
\input xy
\xyoption{all}

\newcommand{\mD}{\mathbb D}
\newcommand{\mR}{\mathbb R}
\newcommand{\mZ}{\mathbb Z}
\newcommand{\mH}{\mathbb H}
\newcommand{\mI}{\mathbb I}

\title{Rigidity result for certain 3-dimensional singular spaces and
their fundamental groups.}

\author{Jean-Fran\c{c}ois Lafont}

\theoremstyle{definition}
\newtheorem{Def}{Definition}[section]

\theoremstyle{proposition}
\newtheorem{Lem}{Lemma}[section]
\newtheorem{Prop}{Proposition}[section]

\newtheorem{Claim2}{Claim}

\theoremstyle{plain}

\newtheorem{Thm}{Theorem}[section]
\newtheorem{Cor}{Corollary}[section]

\theoremstyle{remark}

\newtheorem*{Prf}{Proof}

\begin{document}

\maketitle

\begin{abstract}
In this paper, we introduce a particularly
nice family of $CAT(-1)$ spaces, which we call hyperbolic P-manifolds.  
For $X^3$ a simple, thick hyperbolic 
P-manifold of dimension $3$, we show that certain subsets of the 
boundary at infinity of the universal cover of $X^3$ are characterized
topologically.  Straightforward consequences include a version of
Mostow rigidity, as well as 
quasi-isometry rigidity for these spaces.
\end{abstract}

\section{Introduction.}

In this paper, we prove various rigidity results for certain
3-dimensional hyperbolic P-manifolds.  These spaces form a family of 
stratified metric spaces built up from hyperbolic manifolds 
with boundary (a precise definition is given below).
The main technical tool is an analysis of the boundary at infinity
of the spaces we are interested in.  We introduce the notion of a
{\it branching} point in an arbitrary topological space, and show how 
branching points in the boundary at infinity can be used to determine
both the various strata and how they are pieced together.  The argument
for this relies on a result which might be of independant interest: a
version of the Jordan separation theorem that applies to maps $S^1
\rightarrow S^2$ which are not necessarily injective.  In this 
section, we introduce the objects we are interested in, 
provide some basic definitions, and state the theorems 
we obtain.  The proof of the main theorem will be given in 
section 2.  The various applications will be discussed in section 3.
We will close the paper with a few concluding remarks and open
questions in section 4.

\vskip 10pt

\centerline{\bf Acknowledgments.}

\vskip 5pt

The results contained in this paper were part of the author's thesis,
completed at the University of Michigan, under the guidance of professor 
R. Spatzier.  The author would like to thank his advisor for his help 
throughout the author's graduate years.  The author would also like to 
thank professor J. Heinonen for taking the time
to proofread his thesis and to comment on the results contained therein.
Finally, the author gratefully acknowledges the referee's efforts in 
catching numerous minor errors in the first draft of this paper.

\subsection{Hyperbolic P-manifolds.}

\begin{Def}
We define a closed $n$-dimensional {\it piecewise manifold}
(henceforth abbreviated to P-manifold) to be a topological
space which has a natural stratification into pieces which are
manifolds. More precisely, we define a $1$-dimensional P-manifold
to be a finite graph.
An $n$-dimensional P-manifold ($n\geq 2$) is defined inductively as a 
closed pair
$X_{n-1}\subset X_n$ satisfying the following conditions:

\begin{itemize}
\item Each connected component of $X_{n-1}$ is either an
$(n-1)$-dimensional P-manifold, or an $(n-1)$-dimensional manifold.
\item The closure of each connected component of $X_n-X_{n-1}$
is homeomorphic to a compact orientable $n$-manifold with
boundary, and the homeomorphism takes the component of $X_n-X_{n-1}$ to the interior of the $n$-manifold; the closure of such a component will be called a {\it chamber}.
\end{itemize}

\noindent Denoting closure of the connected components of $X_n-X_{n-1}$ by $W_i$, we observe that we have a natural map $\rho: \coprod
\partial W_I \longrightarrow X_{n-1}$ from the disjoint union
of the boundary components of the chambers to the subspace
$X_{n-1}$.  We also require this map to be surjective, and a
homeomorphism when restricted to each component. The P-manifold is
said to be \emph{thick} provided that each point in $X_{n-1}$ has
at least three pre-images under $\rho$. We will
henceforth use a superscript $X^n$ to refer to an $n$-dimensional
P-manifold, and will reserve the use of subscripts $X_{n-1}\ldots
X_1$ to refer to the lower dimensional strata.  For a thick
$n$-dimensional P-manifold, we will call the $X_{n-1}$ strata the
{\it branching locus} of the P-manifold.
\end{Def}

Intuitively, we can think of P-manifolds as being ``built'' by
gluing manifolds with boundary together along lower dimensional
pieces. Examples of P-manifolds include finite graphs and soap 
bubble clusters.  Observe that compact manifolds can also be 
viewed as (non-thick) P-manifolds.
Less trivial examples can be constructed more or less arbitrarily
by finding families of manifolds with homeomorphic boundary and
glueing them together along the boundary using arbitrary
homeomorphisms. We now define the family of metrics we are
interested in.

\begin{Def}
A Riemannian metric on a 1-dimensional P-manifold (finite graph)
is merely a length function on the edge set.  A Riemannian metric
on an $n$-dimensional P-manifold $X^n$ is obtained by first
building a Riemannian metric on the $X_{n-1}$ subspace, then
picking, for each chamber $W_i$ a Riemannian metric with totally
geodesic boundary satisfying that the gluing map $\rho$ is an isometry. 
We say that a Riemannian metric on a P-manifold is
hyperbolic if at each step, the metric
on each $W_i$ is hyperbolic.  
\end{Def}

\begin{figure}
\begin{center}
\includegraphics[width=2in]{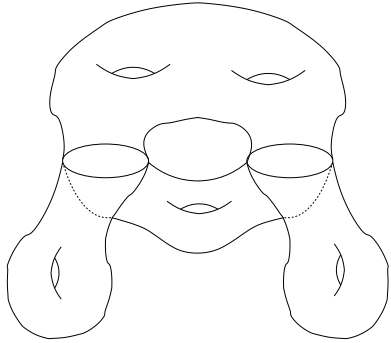}
\caption{Example of a simple, thick P-manifold.}
\end{center}
\end{figure}

A hyperbolic P-manifold $X^n$ is automatically a locally $CAT(-1)$ 
space (see Chapter II.11 in Bridson-Haefliger [5].  Furthermore, the lower dimensional 
strata $X_i$ are all totally geodesic subspaces of $X^n$.  In particular,
the universal cover $\tilde X^n$ of a hyperbolic P-manifold $X^n$ is a $CAT(-1)$ space (so is automatically $\delta$-hyperbolic), and 
has a well-defined boundary at infinity $\partial
^\infty \tilde X^n$. Finally we note that the fundamental group $\pi_1(X_n)$ is a $\delta$-hyperbolic group.

We also note that examples of hyperbolic P-manifolds are easy to
obtain.  In dimension two, for instance, one can take multiple
copies of a compact hyperbolic manifold with totally geodesic
boundary, and identify the boundaries together.  In higher
dimension, one can similarly use the arithmetic constructions of
hyperbolic manifolds (see Borel-Harish-Chandra [2])
to find hyperbolic manifolds with isometric
totally geodesic boundaries.  Gluing multiple copies of these
together along their boundaries yield examples of hyperbolic
P-manifolds.  More complicated examples can be constructed by
finding isometric codimension one totally geodesic submanifolds
in {\it distinct} hyperbolic manifolds.  Once again, cutting the
manifolds along the totally geodesic submanifolds yield hyperbolic
manifolds with totally geodesic boundary, which we can glue together
to build hyperbolic P-manifolds
(see the construction of non-arithmetic lattices by
Gromov-Piatetski-Shapiro [9]). 

\begin{Def}
We say that an $n$-dimensional P-manifold $X^n$ is {\it simple} 
provided its codimension two strata is empty.  In other words, the
$(n-1)$-dimensional strata $X_{n-1}$ consists of a disjoint union
of $(n-1)$-dimensional manifolds.
\end{Def}

An illustration of a simple, thick P-manifold is
given in figure 1.  It has four chambers, and two connected
components in the codimension one strata. 
We point out that P-manifolds show up naturally in the setting
of branched coverings. Starting from a hyperbolic manifold $M$, one can
take a totally geodesic
codimension two subspace $N\subset M$, with the property that $N$ bounds a
codimension one totally geodesic subspace $W$.  One could look at a
ramified cover $\bar M_i$ of $M$, of degree $i$, where the ramification is 
over $N$.  This
naturally inherits a (singular) CAT(-1) metric from the metric on $M$ (note
that
Gromov-Thurston [10] have shown that, when $i$ is sufficiently large, it
is possible to smooth this singular metric to a negatively curved Riemannian
metric).  The pre-image $X\subset \bar M_i$ of the 
subspace $W\subset M$ will be a simple hyperbolic P-manifold, isometrically 
embedded in $\bar M_i$ (with respect to the singular metric).  The codimension
one strata in $X$ will consist of a single connected component, isometric
to $N$, and there will be $i$ chambers, each isometric to $W$.

We also point out that related spaces include hyperbolic buildings.  Indeed,
these share the property that they are naturally stratified spaces,
built out of pieces that are isometric
to subsets of hyperbolic space.  The difference lies in that for hyperbolic 
buildings, the boundary of the chambers are not totally geodesic in the 
space.  Hyperbolic buildings have
been studied by Bourdon-Pajot [4]; they obtain quasi-isometric
rigidity for some 2-dimensional hyperbolic buildings (compare to our Theorem
1.3).  Note that, unlike 
hyperbolic P-manifolds, hyperbolic buildings can only exist in low dimensions
($\leq 29$).
This follows from a result of Vinberg [21], showing that compact Coxeter
polyhedra do not exist in hyperbolic spaces of dimension $\geq 30$.  Since
the existence of such polyhedra is a pre-requisite for the existence of
hyperbolic buildings, Vinberg's result immediately implies the desired 
non-existence result.

\vskip 5pt
 
Next we introduce a locally defined topological invariant.
We use $\mD^n$ to denote a closed $n$-dimensional disk, and
$\mD^n_\circ$ to denote its interior.  We also use $\mI$
to denote a closed interval, and $\mI_\circ$ for its interior.

\begin{Def}
Define the 1-{\it tripod} $T$ to be the topological space obtained
by taking the join of a one point set with a three point set.
Denote by $*$ the point in $T$ corresponding to the one point set.
We define the $n$-{\it tripod} ($n\geq 2$) to be the space $T
\times \mD^{n-1}$, and 
call the subset $*\times \mathbb D^{n-1}$ the {\it spine} of the 
tripod $T\times \mathbb D^{n-1}$.  The subset $*\times 
\mathbb D^{n-1}$ separates
$T\times \mathbb D^{n-1}$ into three open sets, which we 
call the {\it open leaves} of the tripod.  The union of an open leaf with the spine will be called a {\it closed leaf} of the spine.
We say that a point $p$ in a topological space $X$ is {\it
$n$-branching} provided there is a topological embedding $f:T\times 
\mathbb D^{n-1}
\longrightarrow X$ such that $p\in f(*\times 
\mathbb D^{n-1}_\circ)$.  
\end{Def}

It is clear that the property of being $n$-branching is invariant
under homeomorphisms.  We show some examples of branching in Figure
2.  Note that, in a simple, thick P-manifold of dimension
$n$, points in the codimension one strata are automatically 
$n$-branching.  One can ask whether this property can be detected
at the level of the boundary at infinity.  This motivates the following:

\begin{figure}
\begin{center}
\includegraphics[width=3.5in]{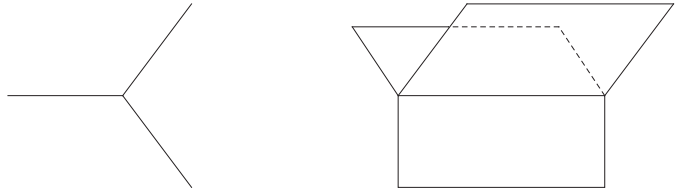}
\caption{Examples of 1-branching and 2-branching.}
\end{center}
\end{figure}

\vskip 5pt

\noindent{\bf Conjecture:} Let $X^n$ be a simple, thick hyperbolic
P-manifold of dimension $n$, and let $p$ be a point in the boundary 
at infinity of
$\tilde X^n$.  Then $p$ is $(n-1)$-branching if and only if $p=\gamma
(\infty)$ for some geodesic ray $\gamma$ contained entirely in a
connected lift of $X_{n-1}$.  

\vskip 5pt

One direction of the above conjecture is easy to prove (see 
Proposition 2.1).  In the case where $n=3$,
we will see that the reverse implication also holds.  Note that in 
general, the (local) topological structure of the boundary at infinity of a $CAT(-1)$ space (or of
a $\delta$-hyperbolic group) is very hard to analyze.  The conjecture
above says that with respect to branching, the boundary of a
simple, thick hyperbolic P-manifold of dimension $n$ is 
particularly easy to understand.

In our proofs, we will make use of a family of nice metrics on the 
boundary at infinity of an arbitrary $n$-dimensional hyperbolic 
P-manifold  (in fact, on the boundary at infinity of any 
CAT(-1) space).

\begin{Def}
Given an $n$-dimensional hyperbolic P-manifold, and a basepoint
$*$ in $\tilde X^n$, we can define a metric on the boundary at
infinity by setting $d_\infty (p,q)=e^{-d(*,\gamma_{pq})}$, where
$\gamma_{pq}$ is the unique geodesic joining the points $p,q$ (and $d$ denotes the distance inside $\tilde X_n$).
\end{Def}

The fact that $d_\infty$ is a metric
on the boundary at infinity of a proper $CAT(-1)$ space follows
from Bourdon (Section 2.5 in [3]).  Note that changing the basepoint from $*$
to $*^\prime$ changes the metric, but that for any $p,q\in
\partial^\infty(X^n)$, we have the inequalities:
$$A^{-1}\cdot d_{\infty, *}(p,q)\leq d_{\infty, *^\prime}(p,q)
\leq A\cdot d_{\infty, *}(p,q)$$ where $A=e^{d(*,*^\prime)}$, and
the subscripts on the $d_\infty$ refers to the choice of basepoint
used in defining the metric.   In particular, different choices
for the basepoint induce the same topology on $\partial ^\infty
\tilde X^n$, and this topology coincides with the standard topology
on $\partial ^\infty \tilde X^n$ (the quotient topology inherited from
the compact-open topology in the definition of $\partial ^\infty \tilde
X^n$ as equivalence classes of geodesic rays in $\tilde X^n$).
This gives us the freedom to select basepoints at our
convenience when looking for {\it topological} properties of the
boundary at infinity.

\subsection{Statement of results, outlines of proofs.}

We will start by proving the following:

\begin{Thm}
Let $X^3$ be a simple, thick hyperbolic P-manifold of dimension $3$.
Then a point $p\in \partial ^\infty \tilde X^3$ is $2$-branching if and only if 
there is a geodesic ray $\gamma _p\subset \tilde X_2$ in a lift of
the $2$-dimensional strata,  
with the property that $\gamma _p(\infty)=p$.
\end{Thm}

This theorem has immediate applications, in that it allows us 
to show several rigidity results for simple, thick hyperbolic P-manifolds
of dimension $3$.  The first application is a version of Mostow
rigidity.

\begin{Thm}[Mostow rigidity.]
Let $X^3_1, X^3_2$ be a pair of simple, thick hyperbolic P-manifolds
of dimension $3$.  Assume
that the fundamental groups of $X^3_1$ and $X^3_2$ are isomorphic.
Then the two P-manifolds are in fact isometric (and the isometry
induces the isomorphism of fundamental groups).
\end{Thm}

This is proved in subsection 3.1, but the idea of the proof is 
fairly straightforward: one uses the
isomorphism of fundamental groups $\Phi:\pi_1(X_1)\longrightarrow
\pi_1(X_2)$ to induce a homeomorphism $\partial ^\infty \Phi:
\partial ^\infty \tilde X_1 \longrightarrow \partial ^\infty \tilde X_2$
between the boundaries at infinity of the respective universal
covers.  Theorem 1.1 implies that the subsets of the
$\partial ^\infty \tilde X_i$ corresponding to the lifts of
the branching loci
are homeomorphically identified.  A separation argument ensures
that the boundaries of the various chambers are likewise identified.
One then uses the dynamics of the $\pi_1(X_i)$ actions on the 
$\partial ^\infty \tilde X_i$ to ensure that corresponding chambers
have the same fundamental groups.  Mostow rigidity for hyperbolic
manifolds with boundary allows us to conclude that, in the quotient,
the corresponding chambers are isometric.  Finally, the dynamics
also allows us to identify the gluings between the various chambers,
yielding the theorem.

As a second application, we consider groups which are
quasi-isometric to the fundamental group of a simple, thick hyperbolic 
P-manifold of dimension $3$.  We obtain:

\begin{Thm}
Let $X$ be a simple, thick hyperbolic P-manifold of dimension $3$.
Let $\Gamma$ be a group quasi-isometric to $\pi_1(X)$.  Then there 
exists a short exact sequence of the form:
$$1\rightarrow F\rightarrow \Gamma \rightarrow \Gamma ^\prime 
\rightarrow 1$$
where $F$ is a finite group, and $\Gamma ^\prime \leq Isom(\tilde X)$
is a subgroup of the isometry group of $\tilde X$ with $\tilde X/
\Gamma ^\prime$ compact.
\end{Thm}

Here our argument relies on showing that any
quasi-isometry of the corresponding P-manifold is in fact a bounded
distance from an isometry.  The key idea is that (by Theorem 1.1) 
any quasi-isometry is bounded distance from a quasi-isometry which
preserves the chambers.  For compact hyperbolic
manifolds {\it with totally geodesic boundary}, a ``folklore theorem''
states that quasi-isometries of the universal cover are a
bounded distance from isometries (note that the corresponding 
statement for {\it closed} hyperbolic manifolds is {\it false});
see the footnote on pg. 648 in Kapovich-Kleiner [12].  
Assuming this result (whose proof we sketch out in section 3.2), 
we then show that the bounded
distance isometries can glue together to give a global isometry
which is still at a bounded distance from the original 
quasi-isometry.
From such a statement, standard methods yield a quasi-isometry
classification.

\section{The main theorem.}

In this section, we provide a proof of Theorem 1.1.  We start
by noting that one direction of the conjecture stated in the 
introduction is easy to prove:

\begin{Prop} 
Let $X^n$ a simple, thick $n$-dimensional 
P-manifold, and let $p$ be a point in the boundary at infinity of
$\tilde X^n$. If $\gamma$ is a geodesic ray contained entirely in a
connected lift $\tilde B$ of $X_{n-1}$, then $\gamma(\infty)$ is 
$(n-1)$-branching.
\end{Prop}

\begin{Prf}
This is easy to show: by the thickness hypothesis, there are at 
least three distinct chambers $\tilde W_i$ containing $\tilde B$
in their closure.  For each of these chambers, we can consider the
various boundary components of $\tilde W_i$.  To each boundary
component distinct from $\tilde B$, we can again use the thickness
hypothesis to find chambers incident to each of the boundary 
components.  Extending this procedure, we see that we can find 
three totally geodesic subset in $\tilde X^n$ glued together 
along the codimension one strata $\tilde B$.  
Furthermore, the simplicity assumption implies that $\tilde B$ is
isometric to $\mH^{n-1}$, while each of the three totally geodesic
subsets is isometric to a ``half'' $\mH^n$.  This implies that, 
in the boundary at infinity, there are three embedded disks 
$\mD^{n-1}$ glued along their boundary to
$S^{n-1}\cong \partial ^\infty \tilde B$.  It is now immediate
that $\gamma (\infty)$ is $(n-1)$-branching.
\end{Prf}

For the reverse implication, we will need a strong form of the
Jordan separation theorem.  The proof of this theorem is the
only place where the condition $n=3$ is used.

\begin{Thm}[Strong Jordan separation]
Let $f:S^1\rightarrow S^2$ be a continuous map, and let $I\subset 
S^1$ be the set of injective points (i.e. points $p\in S^1$ with the
property that $f^{-1}(f(p))=\{p\}$).  If $I$ contains an open set $U$,
and $q\in U$, then:
\begin{itemize}
\item $f(S^1)$ separates $S^2$ into open subsets (we write
$S^{2}-f(S^1)$ as a disjoint union $\amalg U_i$, with each $U_i$ open),
\item there are precisely two open subsets $U_1$, $U_2$ in the 
complement of $f(S^1)$ which contain $p:=f(q)$ in their closure.
\item if $F:\mD^{2}\rightarrow S^{2}$ is an extension of the map $f$ 
to the closed ball, then $F(\mD^{2})$ surjects onto either $U_1$ or $U_2$.
\end{itemize}
\end{Thm}

Before starting with the proof, we note that this theorem clearly 
generalizes the classical Jordan separation theorem in the plane
(corresponding to
the case $I=S^1$).  The author does not know whether the hypotheses on
$I$ can be weakened to just assuming that $I$ is measurable.  

\begin{Prf}
We start out by noting that the map $f(S^1)$ cannot surject onto $S^2$.
Since $I$
is assumed to contain an open set, we can find an $\mI\subset I\subset S^1$,
i.e. a (small) interval on which $f$ is injective.
Since $f$ is injective on $\mI$, one can find a small 
closed ball in $S^2$ with the property that the ball intersects $f(S^1)$
in a subset of $f(\mI_\circ)$.   But an imbedding of a 1-dimensional
space into a 2-dimensional space has an image which must have zero 
measure, so in particular, there is a point in the closed ball that is
not in the image of $f(\mI_\circ)$.

Since $f$ is not surjective, we use stereographic projection to view $f$
as a map into $\mR^2$.  A well known theorem (which is a consequence of the $2$-dimensional Schoenflies theorem) now tells 
us that, given any embedded arc in the plane, there is a homeomorphism
of the plane taking the arc to a subinterval of the $x$-axis (this follows 
for instance from Theorem III.6.B in Bing [1]).  Applying this 
homeomorphism we can assume that $f$ maps the interval $\mI$ to the 
$x$-axis.  Now let $x_1,x_2$ be a pair of points lying slightly above and
slightly below the image $f(\mI)$.  If the points $x_i$ are close enough
to the $x$-axis, we can find a path $\eta$ which 
intersects the $f(\mI)$ transversaly in a single point, joins $x_1$ to 
$x_2$, and has no
other intersection with $f(S^1)$.  Now perturb the map $f$, away from 
$f(\mI)$, so that it is PL.  If the perturbation is slight enough, 
the new map $g$ will be homotopic to $f$ in the complement of the $x_i$.
Furthermore, $\eta$ will intersect the map $g$ transversaly in 
precisely one point.  It is now classical that the map $g$ must represent
distinct elements in $H_1(\mR^{2}-x_1)\cong \mZ$ and $H_1(\mR^{2}-x_2)\cong \mZ$ 
(and in fact, that the integers it represents differ by one).  Since $g$
is homotopic to $f$ in the complement of the $x_i$, the same holds for
the map $f$.  In 
particular, the connected components in $S^2-f(S^1)$
containing $x_1$ and $x_2$ are 
distinct, giving the first two claims.  Furthermore, $f$ represents 
a non-zero
class in one of the $H_1(\mR^{2}-x_i)$, giving us the third claim.
\end{Prf}

Note that geodesic rays $\gamma$ which are not asymptotic
to a ray contained in a lift of the branching locus are of one
of two types:
\begin{itemize}
\item either $\gamma$ eventually stays trapped in a $\tilde W_i$,
and is not asymptotic to any boundary component, or
\item $\gamma$ passes through infinitely many connected lifts
$\tilde W_i$.
\end{itemize}
In the next proposition, we deal with the first of these
two cases.  Let us first introduce some notation.  Given a point $x\in
\tilde X^3$, we denote by $\pi_x:\partial ^\infty \tilde X^3
\longrightarrow lk(x)$ the geodesic projection from the boundary onto
the link at the point $x$.    Recall that the {\it link} of a point
in a piece-wise hyperbolic CAT(-1) space is a small
metric sphere of radius $\epsilon$ centered at the point.  If 
$\epsilon$ is small enough, the link is unique upto homeomorphism.

We denote by $I_x$ the set $\{p\in
lk(x) \hskip 5pt : \hskip 5pt |\pi_x^{-1}(p)|=1\}\subset
lk(x)$, in other words, the set of points in the link where
the projection map is actually injective. The importance of this
set lies in that it consists of those directions (points in the
boundary) where injectivity can be detected {\it from the point
$x$}.

\begin{Prop}
Let $X^3$ be a simple, thick 3-dimensional hyperbolic P-manifold.
Let $\gamma\subset \tilde X^3$ be a geodesic ray lying entirely
within a connected lift $\tilde W$ of a chamber $W$, and not
asymptotic to any boundary component of $\tilde W$.  Then
$\gamma (\infty)$ is {\bf not} 2-branching.
\end{Prop}

\begin{Prf}

We start by observing 
that, by our hypothesis, we can take any $x\in \gamma$ as a
basepoint, and $\pi_x(\gamma(\infty))$ will lie within $I_x$
(since by hypothesis $\gamma$ lies entirely within $\tilde W_i$).
Now assume, by way of contradiction, that $\gamma(\infty)\in 
\partial ^\infty \tilde X^3$ is 2-branching.  Then we have an 
injective map
$f:T\times \mI \longrightarrow \partial ^\infty \tilde X^3$
such that $\gamma(\infty)\in f(*\times \mI_\circ)$. Consider
the composite map $\pi_x\circ f:T\times \mI\longrightarrow
lk(x)$ into the link at $x$.
Since $x$ lies in a chamber,
we have $lk(x)\cong S^2$. Now note that the composite
map $\pi_x\circ f$ must be injective on the set $(\pi_x\circ
f)^{-1}(I_x)$.  Indeed, by the definition of $I_x$, those are the
points $p$ in $lk(x)$ which have a unique pre-image under
$\pi_x$.  Hence, if the composite map $\pi_x\circ f$ has {\it
more} than one pre-image at such a point $p$ , it would force the
map $f$ to have two distinct pre-images at $\pi_x^{-1}(p)\in
\partial ^\infty \tilde X^3$, which violates our assumption that
$f$ is injective.

In order to get a contradiction, we plan on showing that the
composite map fails to be injective at some point in 
the set $I_x$.  We start with a few observations on the 
structure of the set $I_x$.

\begin{Claim2}
The complement of $I_x$ has the following properties:
\begin{itemize}
\item it consists of a countable union of {\it open} disks 
$U_i$ in $S^2$,
\item the $U_i$ are the interiors of a family of pairwise 
disjoint closed disks,
\item the $U_i$ are dense in $S^2$.
\end{itemize}
\end{Claim2}

\begin{Prf}
If $\pi_x$ fails to be injective at a point $p\in lk(x)$,
then there are two distinct geodesic rays emanating from $x$, in
the direction $p$.  Since $x$ lies within a chamber these two
geodesic rays must agree up until some point, and then diverge.
This forces these geodesic rays to intersect the branching locus.

This immediately tells us that the set $I_x$ is the projection of
$\partial ^\infty \tilde W$ onto the link.  Note that this is a 
homeomorphism, and since $\partial ^\infty \tilde W$ is a 
Sierpinski carpet, we immediately get all three claims.
\end{Prf}

Since the point $\pi_x(\gamma (\infty))$ lies in $I_x-\cup (\partial
U_i)$, we would like to get some further information about the 
density of the $U_i$ away from the set $\cup (\partial U_i)$.  

\begin{Claim2}
For any point
$p\in I_x -\cup (\partial U_i)$, and any neighborhood $N_p$ of $p$, 
there exist arbitrarily small $U_i$ with $U_i\subset N_p$.
\end{Claim2}

\begin{Prf}
By density of the $U_i$, we have that for any point $p\in I_x
- \cup (\partial U_i)$, 
arbitrarily small neighborhoods of $p$ must intersect an open disk. To
see that arbitrary small neighborhoods actually {\it contain} an
open $U_i$, we consider the standard measure $\mu$ on the sphere
(identified with the link).  Note that since the measure of the
sphere is finite, for any $\epsilon>0$ there are at most finitely
many $U_i$ with $\mu(U_i)>\epsilon$.  Since the union of the 
boundaries of these $U_i$ form a closed subset of $S^2$, and 
this subset does not contain $p$ (since we assumed $p\in I_x 
-\cup (\partial U_i)$), we have that the distance from $p$ to
the boundaries of these $U_i$ is positive.

In particular, for an
arbitrary neighborhood $N_p$ of $p$, and an arbitrary
$\epsilon>0$, we can find a smaller neighborhood
$N_p^\prime\subset N_p$ with the property that any $U_i$
intersecting $N_p^\prime$ satisfies $\mu(U_i)<\epsilon$. However,
since the $U_i$ are actually {\it round} disks in $S^2$, we
have that $diam(U_i)<C\cdot \mu(U_i)^{\frac {1}{2}}$ (for some
uniform constant $C$), which gives us control of $diam(U_i)$ in
terms of $\mu(U_i)$. So in particular, picking $N_p^\prime$ much
smaller than $N_p$, we can force $diam(U_i)$ to be much smaller
than the distance from $N_p^\prime$ to the boundary of $N_p$.
Hence $U_i\subset N_p$, completing the claim.
\end{Prf}

\begin{Claim2}
The image $(\pi_x\circ f)(\partial(T\times \mI))$ is a bounded
distance away from $\pi_x(\gamma(\infty))$.
\end{Claim2}

\begin{Prf}
First of all, observe that the boundary $\partial (T\times
\mI)$ of the set $T\times \mI$ is compact, forcing $(\pi_x
\circ f)(\partial(T\times \mI))$ to be compact.  Since $f$ is
injective by hypothesis, we must have $\gamma(\infty) \notin
f(\partial (T\times \mI))$, yielding $\pi_x (\gamma(\infty)) \notin
(\pi_x \circ f)(\partial(T\times \mI))$ as we know that $\pi_x(\gamma (\infty))$ lies in the injectivity set $I_x$ . Hence the minimal
distance between $\pi_x (\gamma(\infty))$ and $(\pi_x \circ
f)(\partial(T\times \mI))$ is positive.
\end{Prf}

Now recall that we need to find a point in $I_x$ where the 
composite map $\pi_x\circ f:T\times \mI\rightarrow S^2$ fails to 
be injective.  To do this, we start by observing the following:

\begin{Claim2}
The image of the spine is entirely contained in $I_x$ (i.e.
$(\pi_x\circ f)(*\times \mI_\circ)\subset I_x$).
\end{Claim2}

Heuristically, the idea is that if the claim was false, one would
find a $U_i$ intersecting the image of the spine.  The pre-image
of the boundary of this $U_i$ would look like a tripod $T$ within
the space $T\times \mI$ (see Figure 3).  But $\partial U_i$ lies
within the set $I_x$, so its pre-image should be homeomorphic to a subset
of $S^1$. 

\begin{Prf}
We argue by contradiction.  If not, then there exists a point $q\in
(\pi_x\circ f)(*\times \mI_\circ)$ with the property that $q\notin I_x$.
This implies that $q$ lies in one of the open disks $U_i$.  Note that
we already have a point $p:=\pi_x (\gamma(\infty))$ whose image lies in $I_x -\cup (\partial U_i)$.

\begin{figure}
\begin{center}
\includegraphics[width=3.5in]{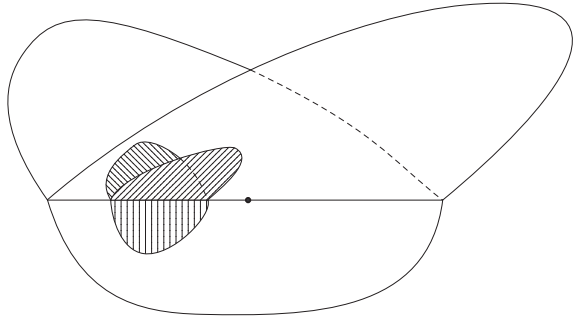}
\caption{Pre-image of a $U_i$ which intersects the spine.}
\end{center}
\end{figure}

Now consider the pre-image $K$ of $\partial U_i$.  Let $L_j$ 
($1\leq j\leq 3$) be the three closed leaves of the tripod, and consider the
intersection $K^\prime := K\cap (L_1\cup L_2)$.  The set $K^\prime$
is the pre-image of $\partial U_i$ for the restriction of the map to
the union $L_1\cup L_2$, hence must separate $p$ and $q$.  $K^\prime$
is a closed subset of $L_1\cup L_2\cong \mD^2$, and since $\partial U_i
\subset I_x$, $K^\prime$ must be homeomorphic to a closed subset of
$\partial U_i\cong S^1$.  This implies that $K^\prime$ consists of 
either a union of intervals, or of a single $S^1$.

We first note that $K^\prime$ {\it cannot} be an $S^1$, for then 
$K^\prime$ would have to equal $K$ (since the map is injective on 
$\partial U_i$).  One could then take a path in the third leaf joining
$p$ to $q$, contradicting the fact that $K$ separates.

So we are left with dealing with the case where $K^\prime$ is a 
union of intervals.  Now let $J\subset K^\prime$ be a subinterval
that separates $p$ from $q$.  Note that such an 
interval must exist, else $K^\prime$ itself would fail to separate.  
Now $J$ not only separates, but also locally separates $T\times \mI$.
Furthermore, $J$ cannot be contained entirely in the spine, so  
restricting and reparametrizing if need be, we can assume that 
there is a subinterval $J_1$ having 
the property that $J_1([0,y])\subset L_1$, and $J_1([x,1])\subset L_2$,
where we are now viewing $J_1$ as a map from $\mI$ into $L_1\cup L_2$,
and $0<x\leq y<1$ (so in particular, $J_1([x,y])$ lies in the spine).  
Observe that since $J_1$ locally separates, we have that 
near $J_1 (x)$, $L_3$ must map into one component of
$S^2-\partial U_i$, while near $J_1(y)$ it must map into the other
component of $S^2-\partial U_i$.  This implies that there is a subinterval
$J_2$ of $K$, lying in $L_3$, and separating the points near $J_1(x)$ 
from those near $J_1(y)$.  But the union $J_1 \cup J_2$ is now
a subset of $K$ homeomorphic to a tripod $T$.
Since $K$ is homeomorphic to a subset of $S^1$, this gives us a contradiction,
completing the claim.  
\end{Prf}

We now focus on the restriction $F_i$ ($1\leq i\leq 3$) of the 
composite map to each of the three closed leafs.
Each $F_i$ is a map from $\mD^2$ to $S^2$, and all three
maps coincide on an interval $\mI\subset S^1=\partial \mD^2$ 
(corresponding to the spine $*\times \mI$).  From Claim 4,
each of the maps $F_i$ is injective on $\mI$. 

\begin{Claim2}
There is a connected open set $W\subset S^2$ with the property
that:
\begin{itemize}
\item at least two of the maps $F_i$ surject onto $W$
\item the closure of $W$ contains the point $\pi_x(\gamma (\infty))$
\end{itemize}
\end{Claim2}

\begin{Prf}
To show this claim, we invoke the strong form of Jordan separation
(Proposition 2.1).  Denote by $G_i$ the restriction of the
map $F_i$ to the boundary of each leaf.  From the strong Jordan 
separation, each $G_i(S^1)$ separates $S^2$, and there are precisely
two connected open sets $U_i,V_i\subset S^2-G_i(S^1)$ which contain 
$G_i(\mI)$ in their closure.  Furthermore, each of the maps $F_i$ 
surjects onto either $V_i$ or $U_i$.

Now if $r$ is small enough, we will have that the ball $D$
of radius
$r$ centered at $\pi_x(\gamma (\infty))$ only intersects $G_i(\mI)$
(this follows from claim 3).  In particular, each path connected
component of $D-G_i(\mI)$ is contained in either $U_i$, or in $V_i$.
Furthermore, by an argument identical to that in Proposition 2.1,
there will be precisely two path connected components $U$,$V$,
of $D-G_i(\mI)$ containing $\pi_x(\gamma (\infty))$ in their closure.
Note that since the maps $G_i$ all coincide on $\mI$, we must have 
(upto relabelling) $U\subset U_i$ and $V\subset V_i$ {\it for each $i$}.
From the strong Jordan separation, we know that each extension $F_i$
surjects onto either $U_i$ or $V_i$, which implies that either $U$
or $V$ lies in the image of two of the $F_i$.  This yields our claim.
\end{Prf}

\begin{Claim2}
Let $V\subset S^2$ be a connected open set, containing the 
point $\pi_x(\gamma (\infty))$ in it's closure.  Then $V$ contains
a connected open set $U_j$ lying in the complement of the set $I_x$.
\end{Claim2}

\begin{Prf}
We first claim that the connected open set $V$ contains a point from 
$I_x$.  Indeed, if not, then $V$ would lie entirely in the complement
of $I_x$, hence would lie in some $U_i$.  Since $\pi_x(\gamma (\infty))$
lies in the closure of $V$, it would also lie in the closure of $U_i$,
contradicting the fact that $\gamma$ is {\it not} asymptotic to any
of the boundary components of the chamber containing $\gamma$.

So not only does $V$ contain the point $\pi_x(\gamma (\infty))$ in its
closure, it also contains some point $q$ in $I_x$.  We claim it in fact 
contains a point in $I_x-\cup (\partial U_i)$.  If $q$ itself lies in
$I_x-\cup (\partial U_i)$ then we are done.  The other possibility is 
that $q$ lies in the boundary of one of the $U_i$.  Now
since $V$ is connected, there exists a path $\eta$ joining $q$ to 
$\pi_x(\gamma (\infty))$.  Now assume that $\eta \cap \big(I_x-\cup 
(\partial U_i) \big)
=\{\pi_x(\gamma (\infty))\}$.  Let $\bar U_i$ denote the closed disks
(closure of the $U_i$),
and note that the complement of the set $I_x-\cup (\partial U_i)$ is the set $\cup (\bar U_i)$.

A result of Sierpinski [20] states the following: let
$X$ be an arbitrary topological space, $\{C_i\}$ a countable
collection of disjoint path connected closed subsets in $X$.  Then
the path connected components of $\cup C_i$ are precisely the
individual $C_i$. 

Applying Sierpinski's result to the set $\cup (\bar U_i)$ shows 
that the path connected component of this union are precisely the
individual $\bar U_i$.  So if $\eta \cap \big(I_x-\cup 
(\partial U_i) \big)=\{\pi_x(\gamma (\infty)\}$, we see that the
path $\eta$ must lie entirely within the $\bar U_i$ containing
$q$.  This again contradicts the fact that $\pi_x(\gamma (\infty))\notin
\cup (\partial U_i)$.
Finally, the fact that $V$ contains a point in $I_x-\cup (\partial U_i)$
allows us to invoke Claim 2, which tells us that there is some $U_j$ which 
is contained entirely within the set $V$, completing our argument.
\end{Prf}

Finally, we note that Claim 6 shows that we must have one of the
connected open components $U_j$ of $S^2-I_x$ lying in the image of
at least two distinct leaves.  In particular, the boundary of the 
set $U_i$ is an $S^1$ which lies in the image of two 
distinct closed leaves.  Since the boundary lies in the set of injectivity,
$I_x$, the only way this is possible is if the {\it spine} maps to
the boundary of $U_i$.  As the map is injective on the spine, this
implies that the spine $*\times \mI$ contains an embedded copy of $S^1$.
This gives us our contradiction, completing the proof of the 
proposition.
\end{Prf}

We now have to deal with the second possibility: that of geodesic
rays that pass through infinitely many connected lifts $\tilde W_i$.
We start by proving a few lemmas concerning separability properties
for the $\tilde B_i$ and $\tilde W_i$, which will also be usefull for 
our applications.

\begin{Lem}
Let $\tilde B_i$ be a connected lift of the branching locus, and
let $\tilde W_j$, $\tilde W_k$ be two lifts of chambers which are
both incident to $\tilde B_i$.  Then $\tilde W_j- \tilde
B_i $ and $\tilde W_k- \tilde B_i $
lie in different connected components of $\tilde X^3 - \tilde
B_i$.
\end{Lem}

\begin{Prf}
We start with a trivial observation, which will be crucial in the
proof of both this lemma and the following one.  Take any cyclic
sequence $\tilde W_0, \tilde B_0, \tilde W_1, \tilde B_1, \ldots
,\tilde W_r, \tilde B_r$ of distinct connected lifts of chambers
and branching locus with the property that each term is incident
to the following one. Then the union of all these sets forms a
totally geodesic subset of $\tilde X^3$. Furthermore, by a
simple application of Seifert-Van Kampen, we find that this
totally geodesic subset of a simply connected non-positively
curved space has $\pi_1\cong \mathbb{Z}$.  But this is impossible,
so no such sequence can exist.

Now, assume that we have two lifts of chambers $\tilde W_j$,
$\tilde W_k$ which are both incident to a connected lift $\tilde
B_i$, but which lie in the same connected component of $\tilde X^3
- \tilde B_i$.  Then taking a geodesic joining a point in  $\tilde
W_j$ to a point in  $\tilde W_k$ but not intersecting $\tilde
B_i$, we can consider the sequence of (connected lifts of)
chambers and branching locus that the geodesic passes through to
get a sequence as above.  But as we explained, this gives us a 
contradiction.
\end{Prf}

\begin{Lem}
Let $\tilde W_i$ be a connected lift of a chamber, and let $\tilde
B_j$, $\tilde B_k$ be two connected lifts of the branching locus
which are both incident to $\tilde W_i$.  Then $\tilde B_j$ and
$\tilde B_k$ lie in different connected components of $\tilde X^3
- Int(\tilde W_I)$.
\end{Lem}

\begin{Prf}
This proof is identical to the previous one: just interchange the
roles of the connected lifts of chambers and the connected lifts
of the branching locus.
\end{Prf}

Next, we note that, in the setting we are considering, we can push 
the separability properties out to infinity, obtaining that the 
corresponding boundary points separate.

\begin{Lem}
Let $\partial^\infty \tilde B_i$ be the boundary at infinity of a
connected lift of the branching locus, and let $\partial^\infty
\tilde W_j$, $\partial^\infty \tilde W_k$ be the boundaries at
infinity of two lifts of chambers which are both incident to
$\tilde B_i$. Then $\partial^\infty \tilde W_j- \partial^\infty \tilde B_i $ and
$\partial^\infty \tilde W_k- \partial^\infty \tilde B_i$ lie in different connected components
of $\partial^\infty \tilde X^3 - \partial^\infty \tilde B_i$.
\end{Lem}

\begin{Prf}
Let $\eta:[0,1]\longrightarrow \partial^\infty \tilde X^3$ be a
path in the boundary at infinity joining a point in
$\partial^\infty \tilde W_j- \partial^\infty \tilde B_i$ to a point in $\partial^\infty \tilde
W_k - \partial^\infty \tilde B_I$, which avoids $\partial^\infty \tilde B_i$. Fix a basepoint
$p\in \tilde B_i$, and consider the pair of geodesics
$\gamma_i:[0,\infty)\longrightarrow \tilde X^3$ ($i=0,1$)
satisfying $\gamma_i(0)=p$, and $\gamma_i(\infty)=\eta(i)$.

Now observe that, by assumption, $\eta([0,1])\cap \partial^\infty
\tilde B_i=\emptyset$, and as they are both compact subsets, this
forces $d_\infty(\eta([0,1]), \partial^\infty \tilde
B_i)>\epsilon>0$.  So let us consider a covering of $\eta([0,1])$
by open balls of radius $r=\epsilon/4$ in the compactification
$\tilde X^3\cup \partial^\infty \tilde X^3$. Note that these open
balls are all path-connected.  By compactness of $\eta$,
we can extract a finite
subcover $\{U_i\}$ which still covers $\eta$. The union of these
open sets form a neighborhood of $\eta$ in the compactification,
which, by our choice of $r$ cannot intersect $\tilde B_i$.
Furthermore, this neighborhood is connected, and for sufficiently
large $t$, both $\gamma_0(t)$ and $\gamma_1(t)$ lie in the
neighborhood.  By concatenation of paths, we can obtain a path
$\gamma$ which completely avoids $\tilde B_i$, but joins a point
in $\tilde W_i$ to a point in $\tilde W_j$.  However, we have
already shown that the latter two subsets lie in distinct path
components of $\tilde X^3 -\tilde B_i$.  Our claim follows.
\end{Prf}

\begin{Lem}
Let $\partial^\infty \tilde W_i$ be the boundary at infinity
corresponding to a connected lift of a chamber, and let
$\partial^\infty \tilde B_j$, $\partial^\infty \tilde B_k$ be the
boundary at infinity of two connected lifts of the branching locus
which are both incident to $\tilde W_i$. Then $\partial^\infty
\tilde B_j$ and $\partial^\infty \tilde B_k$ lie in different
connected components of $\partial^\infty \tilde X^3 -
(\partial^\infty \tilde W_i - \cup \partial^\infty \tilde B_l)$, where the union is over all $\tilde B_l$ which are boundary components of $\tilde W_i$.
\end{Lem}

\begin{Prf}
Let us start by noting that all of the sets $\partial^\infty
\tilde B_j$ are closed subsets of $\partial^\infty \tilde X^3$.
Let us focus on those $\tilde B_l$ which are the boundary of our
$\tilde W_i$.  By our previous result, each of those separates
within $\tilde X^3$.  So for each of them, we can consider the
union of the components which {\it do not} contain $\tilde W_i$.
Together with the corresponding $\tilde B_l$, these will form a
countable family of closed totally geodesic subsets $C_l$ indexed
by the boundary components of $\tilde W_i$. Consider the
corresponding subsets $\partial ^\infty C_l$ in $\partial ^\infty
\tilde X^3$. Since each of these $C_l$ is totally geodesic, the
corresponding subset $\partial ^\infty C_l$ is a closed subset of
$\partial ^\infty \tilde X^3$.  Furthermore, their union is the
whole of $\partial ^\infty \tilde X^3-(\partial ^\infty \tilde
W_i I - \cup \partial^\infty \tilde B_l)$.  We now claim that the sets $\partial ^\infty C_l$ are
pairwise disjoint.  But this is clear: by construction, we have
that the $\partial^\infty \tilde B_l$ separate $\partial ^\infty
C_l$ from all the other $\partial ^\infty C_k$.  So the distance
from any point $p\in \partial ^\infty C_l$ to any point $q\in
\partial ^\infty C_k$ is at least as large as the distance between
the corresponding $\partial^\infty \tilde B_l$ and
$\partial^\infty \tilde B_k$.  But since the two totally geodesic
subsets $\tilde B_l$ and $\tilde B_k$ diverge exponentially, the
sets $\partial ^\infty \tilde B_l$ and $\partial ^\infty \tilde B_k$
are some positive distance apart.

Finally, let us assume there is some path $\eta:\mI\longrightarrow
\partial^\infty \tilde X^3$ satisfying $\eta(0)\in \partial^\infty
\tilde B_j$, $\eta(1)\in \partial^\infty \tilde B_k$ ($j\neq k$),
and $\eta\cap \partial^\infty \tilde W_i=\emptyset$.  Then $\eta$
is a continuous map that lies entirely in the complement of
$\partial^\infty \tilde W_i$.  Consider the pre-image of the
various closed sets $\partial^\infty \tilde C_r$ under $\eta$.
This provides a covering of the unit interval by a countable
family of disjoint closed sets.  Applying the result of Sierpinski [20]
which was stated in the proof of Proposition 2.2 (Claim 6),
this is impossible unless the covering is by a single set, 
consisting of a single interval.  This concludes our argument.
\end{Prf}

Note that the previous two lemmas allow us to identify separability
properties within the space with separability properties on the
boundary at infinity.  In particular, 
we can talk about a point within the space lying in a
different component from a point at infinity (i.e. the unique
geodesic joining the pair of points intersects the totally
geodesic separating subset, whether this is a $\tilde B_i$ or 
a $\tilde W_i$).  We are now ready to deal with the second
case of theorem 1.1:

\begin{Prop}
Let $X^3$ be a simple, thick 3-dimensional hyperbolic P-manifold.
Let $\gamma \subset \tilde X^3$ be a geodesic that passes through
infinitely many connected lifts $\tilde W_i$.  Then
$\gamma (\infty)$ is {\bf not} 2-branching.
\end{Prop}

\begin{Prf}
The approach here consists of reducing to the situation covered
in proposition 2.2.
We start by re-indexing the various consecutive connected lifts
$\tilde W_i$ that $\gamma$ passes through by the integers.  Fix a
basepoint $x\in \tilde W_0$ interior to the connected lift $\tilde
W_0$, and lying on $\gamma$. Now assume that there is an injective
map $f:T\times \mI\longrightarrow
\partial ^\infty \tilde X^3$ with $\gamma (\infty)\in
f(*\times \mI_\circ)$. 

We start by noting that, between successive connected lifts
$\tilde W_i$ and $\tilde W_{i+1}$ that $\gamma$ passes through,
lies a connected lift of the branching locus, which we denote
$\tilde B_i$.  Observe that distinct connected lifts of the
branching locus stay a uniformly bounded distance apart. Indeed,
any minimal geodesic joining two distinct lifts of the branching
locus must descend to a minimal geodesic in a $W_i$ with endpoints
in the branching locus.  But the length of any such geodesic is
bounded below by half the injectivity radius of $DW_i$, the double of $W_I$ across its boundary.  By setting
$\delta$ to be the infimum, over all the finitely many chambers
$W_i$, of the injectivity radius of the doubles $DW_i$, we have $\delta>0$.
Let $K_i$ be the connected component of $\partial ^\infty \tilde
X^3 -\partial ^\infty \tilde B_i$ containing $\gamma (\infty)$.
Then for every $p\in K_i$ ($i\geq 1$), we have:
$$d_x(p, \gamma (\infty))<e^{-\delta (i-1)}.$$
Indeed, by Lemma 2.1, $\tilde
B_i$ separates $\tilde X^3$ into (at least) two totally geodesic
components. Furthermore, the component containing $\gamma(\infty)$
is {\it distinct} from that containing $x$.  Hence, the distance
from $x$ to the geodesic joining $p$ to $\gamma (\infty)$ is at 
least as large as the
distance from $x$ to $\tilde B_i$.  But the later is bounded below
by $\delta (i-1)$.  Using the definition of the metric at infinity, and picking 
$x$ as our basepoint, our estimate follows.  

Since our estimate shrinks to zero, and since the distance from
$\gamma (\infty)$ to $f(\partial (T\times \mI))$ is positive, we 
must have a
point $q\in f(*\times \mI_\circ)$ satisfying $d_x(q,
\gamma(\infty))>e^{-\delta (i-1)}$
for $i$ sufficiently large. Since $\tilde B_i$ 
separates, we see that for $i$ sufficiently
large, $f(*\times \mI_\circ)$ contains points on both sides of
$\partial ^\infty \tilde B_i$.  This implies that there
is a point $q^\prime \in f(*\times \mI_\circ)$ that lies within 
some $\partial ^\infty
\tilde W_k$. But such a point corresponds to a geodesic ray lying
entirely within $\tilde W_k$, and {\it not} asymptotic to any of
the lifts of the branching locus.  Finally, we note that {\it any}
point in the image $f(*\times \mI_\circ)$ can be considered
$2$-branching, so in particular the point $q^\prime$ is
$2$-branching.  But in the previous proposition, we showed
this is impossible.  Our claim follows.
\end{Prf}

Combining propositions 2.1, 2.2, and 2.3 gives us the result
claimed in theorem 1.1.  We round out this section by making a
simple observation, which will be used in the proofs of 
theorem 1.2 and 1.3.

\begin{Lem}
Let $X^n$ be a simple hyperbolic P-manifold of dimension
at least three
Let $\partial ^\infty \tilde B\subset \partial ^\infty \tilde X^n$
consist of all limit points of geodesics in the branching locus.
If $n\geq 3$, then the maximal path-connected components of
$\partial ^\infty \tilde B$ are precisely the sets of the form
$\partial ^\infty \tilde B_i$, where $\tilde B_i\subset \tilde B$ is a
single connected component of the lifts of the branching locus.
\end{Lem}

\begin{Prf}
Clearly, the sets $\partial ^\infty \tilde B_i$ are closed 
(since the $\tilde B_i$ are totally geodesic) and
path-connected (since each $\tilde B_i$
is an isometrically embedded $\mH^{n-1}$, so the corresponding
$\partial ^\infty \tilde B_i\cong S^{n-2}$).  
Now let $\tilde B_i, \tilde B_j$ be distinct connected
components of $\tilde B$.  We are left with showing that $\partial
^\infty \tilde B_i\cap \partial ^\infty \tilde B_j=\emptyset$.
Consider a geodesic $\gamma$ joining $\tilde{B_i}$ to
$\tilde{B_j}$. Since they are distinct connected lifts of the
branching locus, this geodesic must intersect a $\tilde W_k$. 
By lemma 2.4, a proper subset of $\partial ^\infty \tilde W_k$ 
separates $\partial ^\infty\tilde B_i$ from $\partial
^\infty\tilde B_j$.  In particular, this forces the latter two
sets to be disjoint. To conclude, we apply the result of Sierpinski [20] 
(stated in the proof of Proposition 2.2, Claim 6).
This concludes the proof of the Lemma.
\end{Prf}

\section{Applications: rigidity results.}

\subsection{Mostow rigidity and consequences.}

In this section, we provide a proof of Mostow rigidity for simple,
thick, hyperbolic P-manifolds of dimension $3$ (Theorem 1.2).  We also
mention some immediate consequences of the main theorem.

\begin{Prf}
We are given a pair $X_1$, $X_2$ of simple, thick, hyperbolic 
P-manifolds of dimension $3$, with isomorphic fundamental groups,
and we want to show that the two spaces are isometric.  We start by 
noting that our isomorphism of the fundamental groups is 
a quasi-isometry, so that we get an induced homeomorphism
$\partial ^\infty \Phi:
\partial ^\infty \tilde X_1\longrightarrow \partial ^\infty \tilde
X_2$ between the boundaries at infinity of the two universal
covers $\tilde X_1$ and $\tilde X_2$.
Let $Y_1\subset \partial ^\infty \tilde X_1$, $Y_2\subset
\partial ^\infty \tilde X_2$ be the set of points in the respective 
boundaries at infinity that are $2$-branching.  Note that since 
$\partial ^\infty \Phi$ is a homeomorphism, and since the property
of being $2$-branching is a topological invariant, we must have 
$(\partial ^\infty \Phi )(Y_1)=Y_2$.  Let $B_{1,i}\subset 
\tilde X_1$, $B_{2,i}\subset \tilde X_2$ be the various connected
lifts of the branching locus.

Theorem 1.1 tells us that we have the equalities
$Y_1=\bigcup (\partial ^\infty \tilde B_{1,i})$, 
$Y_2=\bigcup (\partial ^\infty \tilde B_{2,i})$.  
In particular, 
$\partial ^\infty \Phi$ must map each path connected component 
of $Y_1$ to a path connected component of $Y_2$.  This implies
(by Lemma 2.5)
that $\partial ^\infty \Phi$ induces a bijection between the lifts
$\tilde B_{1,i}$ and the lifts $\tilde B_{2,i}$.  Furthermore,
the homeomorphism $\partial ^\infty \Phi$ must map the complement
of the set $Y_1$ to the complement of the set $Y_2$.  Note that
in $\partial ^\infty \tilde X_1$ and $\partial ^\infty \tilde X_2$,
the complements of the sets $Y_1$, $Y_2$ will have path components
of the following two types:
\begin{enumerate}
\item path-isolated points, corresponding to geodesic rays that pass
through infinitely many $\tilde W_i$, and
\item non-path-isolated points, corresponding to
geodesic rays that eventually lie entirely within a fixed $\tilde
W_i$ (and are not asymptotic to a boundary component).
\end{enumerate}

We note that there are uncountably many of the former, but only
countably many path connected components of the latter.  In
particular, our homeomorphism {\it cannot} map a non-isolated
point to an isolated point.  Hence our homeomorphism provides us
with a bijection from the set of connected lifts of chambers in
$\tilde X_1$ to the set of connected lifts of chambers in $\tilde
X_2$. 

The next claim is that if a lift of a chamber $\tilde W_1\subset 
\tilde X_1$ corresponds to a lift of a chamber $\tilde
W_2\subset X_2$, that they are in fact isometric.  To see this,
we consider the chamber $W_1\subset X_1$,$W_2\subset X_2$ whose
lifts we are dealing with, and note that they have isomorphic
fundamental groups.  Indeed, consider the action of the fundamental
groups of the two P-manifolds on their boundary at infinity.  Then
the fundamental group of a chamber $W_i$ can be identified with 
the stabilizer of $\tilde W_i$ for the action of $\pi_1(X_1)$ as 
deck transformations. We would like to identify $\pi_1(W_i)$ from
the boundary at infinity.  This is the content of the following:

\vskip 5pt
\noindent{\bf Assertion:} The stabilizer of the lift of 
a chamber $\tilde W_i$ coincides with the stabilizer of the 
set $\partial ^\infty \tilde W_i$ in the boundary at infinity.
The respective actions are those of $\pi_1(X)$ as deck 
transformations on $\tilde X$, and the corresponding induced 
action on the boundary at infinity.
\vskip 5pt

To see this, we note that the stabilizer of $\tilde W_i$ will 
clearly stabilize $\partial ^\infty \tilde W_i$.  Conversely, 
assume that we have a non-trivial element $\alpha$ in $\pi_1(X)$ 
which stabilizes $\partial ^\infty \tilde W_i$.  Note that, 
the $\tilde B_i$ must be permuted by any isometry, and
from Lemma 2.1 they separate $\tilde X$ into the various 
lifts of chambers.  Hence it is sufficient to exhibit a point 
in $\tilde W_i$ whose image under $\alpha$ is also in 
$\tilde W_i$.

Note that if $\alpha$ stabilizes $\partial ^\infty 
\tilde W_i$, then so do all its powers.  Since $\alpha$ acts 
hyperbolically on the boundary
at infinity, this implies that the sink/source of the $\alpha$
action lies in the set $\partial ^\infty \tilde W_i$.  Hence 
$\alpha$ stabilizes a geodesic $\gamma$
lying entirely in $\tilde W_i$
(joining the sink and source of the $\alpha$ action on the 
boundary at infinity). There are now two possibilities: either
$\gamma$ lies in the interior of $\tilde W_i$ and we are
done, or $\gamma$ lies on the boundary.  If $\gamma$ lies on 
the boundary, then we have that $\alpha$ must stabilize that 
boundary component, call it $\tilde D$.  Now pick a point $q$ 
in $\partial ^\infty \tilde W_i$ which is {\it not} on 
$\partial ^\infty \tilde D$,
and let $\eta$ be a geodesic from a point in $D$ to the point
$q$. Since $\alpha$ stabilizes $D$, and stabilizes $\partial 
^\infty \tilde W_i$ it maps $\eta$ to a geodesic ray emanating
from a point in $D$, and having endpoint {\it not} on 
$\partial ^\infty \tilde D$.  In particular, $\alpha$ maps 
a point in the interior of $\tilde W_i$ (namely an interior 
point on the ray $\eta$) to another interior point.  As we
remarked earlier, this implies that $\alpha$ stabilizes 
$\tilde W_i$, giving us the assertion.

\vskip 5pt

From the assertion, we now have the desired claim that
if $\tilde W_1\subset \tilde X_1$ corresponds to a 
$\tilde W_2\subset \tilde X_2$, then the chambers $W_1$
and $W_2$ have isomorphic fundamental groups.
Mostow rigidity for hyperbolic manifolds with boundary 
(see Frigerio [7]) now allows us to conclude
that the $W_1$ is isometric to $W_2$, and that the isometry 
induces the isomorphism given above.  Lifting this isometry, we
see that there is an isometry of $\tilde W_1$ to $\tilde W_2$
which induces the isomorphism between the two respective
stabilizers.

Next we discuss how the isometries on the lift of the chambers
glue together to give a global isometry.  We first need to 
ensure that adjacent chambers in $X_1$ map to adjacent chambers
in $X_2$.  Note that two chambers in $X_1$ are adjacent if and
only if there is a unique $B_{1,i}$ separating them.  But by
Lemma 2.3, this can be detected on the level of the boundary at 
infinity.  Since the $B_{1,i}$ map bijectively to the $B_{2,i}$,
there will be a unique $B_{1,i}$ separating a pair of chambers
if and only if there is a unique $B_{2,i}$ separating the 
corresponding chambers in $X_2$.  This implies that incident
chambers map to incident chambers.  Finally, we can recognize 
the fundamental group of the common codimension one manifold $B_i$
in terms of the sink/source dynamics of the action of the 
fundamental group of each chamber on the corresponding boundary
component.  This also allows us to recognize the subgroups of the
$\pi_1(W_{1,i})$ and $\pi_1(W_{1,j})$ that get identified.  
Equivariance of the homeomorphism ensures that the corresponding 
image groups get identified in precisely the same way, which 
implies that the corresponding lifts of the chambers are glued
together in an equivariant, isometric manner.
Finally, we see that there is an equivariant isometry between 
the universal covers $\tilde X_1$ and $\tilde X_2$, which gives 
us our desired claim.  It is clear from our construction that the
isometry we obtain induces the original isomorphism between the
fundamental groups.
\end{Prf}

We point out two immediate (and standard) corollaries:

\begin{Cor}
Let $X^3$ be a simple, thick hyperbolic P-manifold of dimension
$3$, $\Gamma$ its fundamental group.  Then the outer automorphism
group $Out(\Gamma)$ is a finite group, isomorphic to $Isom(X^3)$ 
(the isometry group of the P-manifold).
\end{Cor}

\begin{Cor}
Let $X^3$ be a simple, thick hyperbolic P-manifold of dimension
$3$, $\Gamma$ its fundamental group.  Then $\Gamma$ is a  
co-Hopfian group.
\end{Cor}

Concerning corollary 3.1, we remark that Paulin [16] has shown that 
a $\delta$-hyperbolic group with infinite outer automorphism group 
splits over a virtually cyclic group.  As for corollary 3.2, we point
out that Sela has shown that torsion-free $\delta$-hyperbolic 
groups are Hopfian [18], and that a non-elementary torsion-free 
$\delta$-hyperbolic group is co-Hopfian if and only if it is 
freely indecomposable [17].

Now let $\Sigma _3$ consist of those groups 
which arise as the fundamental group of a simple, thick hyperbolic
P-manifold of dimension $3$.  Note that every group in $\Sigma_3$ 
arises as the fundamental group of a graph of groups, induced
by the decomposition of the P-manifold into its chambers (see Serre
[19] for definitions).  
Furthermore, the gluings between the chambers are encoded in
the morphisms attached to each edge in the graph of groups.
A purely group theoretic reformulation of Mostow rigidity
is the following:

\begin{Cor}[Diagram Rigidity]
Let $H_1,H_2$ be groups in $\Sigma _3$.  Then $H_1\cong H_2$
if and only if there is an isomorphism between
the underlying graph of groups with the property that:
\begin{itemize}
\item the isomorphism takes vertex groups to isomorphic vertex groups, 
\item isomorphisms can be chosen between the vertex groups 
which intertwine all the edge morphisms.
\end{itemize}
\end{Cor}

This result essentially asserts that the ``structure'' of the
graph of groups that yield groups in $\Sigma _3$ is in fact unique.  
For related results, we refer to Forester [8] (see also Guirardel [11]).

\subsection{Quasi-isometry rigidity.}

In this section, we provide a proof of Theorem 1.3, giving
a quasi-isometry classification for fundamental groups of 
simple, thick hyperbolic P-manifolds of dimension $3$.
In proving this theorem, we will use the following well known
result (for a proof, see Proposition 3.1 in Farb [6]):

\begin{Lem}
Let $X$ be a proper geodesic metric space, and assume that every
quasi-isometry from $X$ to itself is in fact a bounded distance
from an isometry. Furthermore, assume that a finitely generated
group $G$ is quasi-isometric to $X$.  Then there exists a
cocompact lattice $\Gamma\subset Isom(X)$, and a finite group $F$
which fit into a short exact sequence:
$$0\longrightarrow F\longrightarrow G\longrightarrow \Gamma
\longrightarrow 0$$
\end{Lem}

So to prove the theorem, it is sufficient
to show that any quasi-isometry of a simple, thick P-manifold
of dimension $3$ is a bounded distance away from an
isometry.  In order to do this, we begin by recalling a well
known ``folklore'' result.  Proofs of this have been given at
various times by Farb, Kapovich, Kleiner, Leeb, Schwarz,
Wilkinson, and others, though no published proof exists (both B.
Kleiner and B. Farb were kind enough to e-mail us their arguments,
which we sketch out below).

\begin{Prop}
Let $M^3$ be a compact hyperbolic 3-manifold with
totally geodesic boundary (non-empty).  Then any quasi-isometry of
the universal cover $\tilde M^3$ is a finite distance from an
isometry. 
\end{Prop}

The idea of the argument is to repeatedly reflect $\tilde M^3$
through the totally geodesic boundary components to get a copy of
$\mathbb{H}^3$, tiled by copies of $\tilde M^3$.  Now given a
quasi-isometry of the original $\tilde M^3$, we can extend to a
quasi-isometry of all of $\mathbb{H}^3$, which has the special
property that it preserves the union of the boundaries (as sets).
This quasi-isometry extends to a quasi-conformal homeomorphism of
the boundary $S^2$ that interchanges certain families of $S^1$
(the points at infinity corresponding to the various boundaries).
Using the fact that this homeomorphism preserves a family of
circles containing nested circles of arbitrarily small size, one
shows that the quasi-conformal homeomorphism is in fact conformal.
This implies that there is an isometry of $\mathbb{H}^3$ which is
bounded distance from the original quasi-isometry.  Furthermore,
by construction, this isometry preserves our original $\tilde
M^3$.

Now assuming the preceding folklore theorem, we proceed to
give a proof of Theorem 1.3:

\begin{Prf}
Let us start by showing our first claim: that any quasi-isometry
of the universal cover $\tilde X$ of a simple, thick P-manifold $X$ 
of dimension $3$ lies a finite distance away
from an isometry.  Notice that our quasi-isometry induces a
self-homeomorphism of the boundary at infinity $\partial ^\infty 
\tilde X$.  Once again, Theorem 1.1 implies that the induced map
on the boundary at infinity acts as a permutation on the set 
of boundaries of connected lifts of the branching locus.

In particular, this forces our quasi-isometry to map each
of the branching strata $\tilde B_i$ {\it inside} the P-manifold to
within finite distance of another branching strata, call it
$\tilde B_i^\prime$. Since under a quasi-isometry we have uniform control
of the distance between the images of geodesics and actual
geodesics, we see that there is a uniform upper bound on the
distance between the image of $\tilde B_i$ and the strata $\tilde B_i^\prime$.
As such, we can modify our quasi-isometry by projecting the images
of each $\tilde B_i$ to the corresponding $\tilde B_i^\prime$.  Since this
projection only moves points by a bounded distance, we have that
the new map is still a quasi-isometry, and is bounded distance
from the one we started with.

So we have now reduced to the case where the quasi-isometry maps
each $\tilde B_i$ into the corresponding $\tilde B_i^\prime$.  Since our induced
homeomorphism on the boundary also permutes the boundaries of the
$\tilde W_i$, we can apply the same projection argument to ensure that
our new quasi-isometry actually maps each $\tilde W_i$ strictly into a
corresponding $\tilde W_i^\prime$. Let us denote this new quasi-isometry
by $f$.  Now Proposition 3.1 forces $\tilde W_i\cong
\tilde W_i^\prime$, and the restriction of our quasi-isometry $f$ to
$\tilde W_i$ is a bounded distance from an isometry
$\phi_i:\tilde W_i\longrightarrow \tilde W_i^\prime$.  Furthermore, as in our
proof of Mostow rigidity, a separation argument ensures that
incidence of the chambers $\tilde W_i$,$\tilde W_j$ forces the corresponding
chambers $\tilde W_i^\prime$ and $\tilde W_j^\prime$ to be incident.

We now want to get a global isometry from the isometries on
chambers.  Observe that, for an incident pair of chambers $\tilde W_i$
and $\tilde W_j$, we can consider the branching strata $\tilde W_i\cap \tilde W_j$.
The image of this map under $f$ is $\tilde W_i^\prime\cap \tilde W_j^\prime
\cong \mH^2$.
Furthermore, we have a pair of isometries $\phi_i,\phi_j$ from
$\tilde W_i\cap \tilde W_j$ to $\tilde W_i^\prime\cap \tilde W_j^\prime$, 
each of which is a
finite distance from the map $f$, so in particular, which must be
a finite distance from each other.  Considering the isometry
$g:=\phi_j^{-1}\circ \phi_i:\tilde W_i\cap \tilde W_j\longrightarrow \tilde W_i\cap
\tilde W_j$, we obtain an isometry of $\tilde W_i\cap \tilde W_j$ which is bounded
distance from the identity.  But the only isometry
of $\mH^2$ which is bounded distance from the identity is the
identity itself.  This allows us to conclude that $\phi_i$ and 
$\phi_j$ are exactly
the same isometry when restricted to $\tilde W_i\cap \tilde W_j$, allowing us to
glue them together.  Since this holds for arbitrary incident
chambers, we can combine all the various isometries into a 
globally defined isometry on $\tilde X$.  

We are left with showing that the resulting isometry
is a bounded distance from the original quasi-isometry.  Note
that, for the time being, we only know that on each lift of
a {\it chamber}
$\tilde W_i$ the isometry is bounded distance $D_i$ from an isometry. 
We still need to deal with the possibility that the individual
$D_i$ might be tending to infinity.

This prompts the question: given that a
$(C,K)$-quasi-isometry is a bounded distance from an isometry, can
we obtain a {\it uniform} upper bound on how large this distance
can get?  We need to obtain a uniform bound for quasi-isometries
of the universal cover of compact hyperbolic manifolds with 
(non-empty) totally geodesic boundary.

In order to answer this, we recall that, for an
arbitrary $(C,K)$ quasi-isometry $f$ on a $CAT(\delta)$ space 
$\tilde W$ ($\delta < 0$),
there is a uniform constant $D:=D_{C,K,\delta}$ (depending solely
on the constants $C,K,\delta$) with the following property.  Given
any bi-infinite geodesic $\gamma$, the distance between the image 
$f(\gamma)$
(which is referred to as a {\it quasi-geodesic}) and the 
bi-infinite geodesic
with endpoints $\partial^\infty f(\gamma(\pm \infty))$ is bounded
above by $D$. Naturally, if the quasi-isometry is bounded distance
from an isometry $\phi$, then the latter geodesic is precisely
$\phi(\gamma)$.  

Now let us assume that we are dealing with a space with the
property that every point has a pair of perpendicular bi-infinite
geodesics $\gamma_1, \gamma_2$
intersecting precisely in $p$.  We will abbreviate this property
to saying that a space has (PBIG).  Then for our isometry, we see that
$\phi(p)=\phi(\gamma_1)\cap \phi(\gamma_2)$, while for our
quasi-isometry we only obtain $f(p)\subset f(\gamma_1)\cap
f(\gamma_2)$. Since $f(\gamma_i)$ lies in the $D$ neighborhood of
$\phi(\gamma_i)$, we see that $f(p)$ lies in the intersection of
the $D$-neighborhoods of a pair of intersecting geodesics, which,
since $\phi$ is an isometry, are in fact perpendicular geodesics.
But such a neighborhood has a diameter that is uniformly bounded
by some constant $D^\prime$ which only depends on $D$ (and hence
on $C,K,\delta$).  

Note that the spaces we are interested in are universal covers
of compact hyperbolic manifolds with non-empty boundary, so it is
not clear that the above property holds.  It is easy to see that
every point is contained in a bi-infinite geodesic, but it is 
less clear that one can find two such geodesics which are 
perpendicular.  For the spaces we are considering, we now make the:

\vskip 5pt

\noindent{\bf Assertion:} There exists a constant $D^{\prime \prime}$
with the property that if $\gamma$ is an inextendable geodesic,
then $d(f(\gamma), \phi(\gamma))\leq D^{\prime \prime}$.

\vskip 5pt

Assuming this assertion, it is easy to obtain the upper bound we
desire.  Indeed, every point in $\tilde W$ is the intersection of a 
pair of perpendicular inextendable geodesics $\gamma_1,\gamma_2$
(which either terminate at the boundary, or extend off to infinity).  
The same argument as before shows that $f(p)$ must lie in the 
intersection of the $D^{\prime \prime}$ neighborhoods of 
$\phi (\gamma_1)$ and $\phi (\gamma _2)$, giving uniform control of
$d(f(p), \phi(p))$

To see that the desired assertion is true, we note that we only
have to deal with geodesic segments with both endpoints on boundary
components, or geodesic rays emanating from a boundary component 
(the case of bi-infinite geodesics having been discussed above).
Now note that, since the boundary components have the property
(PBIG) (indeed, they are totally geodesic $\mH ^2$'s), we have that 
for points $q$ on the boundary
components, $d(f(q),\phi(q))\leq D^\prime$.  

If $\gamma$ is a geodesic segment with both endpoints $q_1,q_2$
on boundary components, then we have that $d(f(\gamma), \eta)
\leq D$, where $\eta$ is the geodesic joining $f(q_1)$ to $f(q_2)$.
However, we also have that $d(f(q_i),\phi(q_i))\leq D^\prime$,
so by convexity of the distance function $d(\eta, \nu)\leq D^
\prime$, where $\nu$ is the geodesic joining $\phi(q_1)$ to 
$\phi(q_2)$.  But that geodesic is precisely $\phi(\gamma)$,
so the triangle inequality yields $d(\phi(\gamma),f(\gamma))
\leq D+D^\prime=:D^{\prime \prime}$, giving the desired upper 
bound for this case.
The case of a geodesic ray with endpoint on a boundary component
follows from an identical argument.

We conclude that we have the desired uniform bound, which implies
that our gluing of the `piecewise' isometries is still
a bounded distance from an isometry.

Now a consequence of every quasi-isometry being finite distance
from an isometry is that any group $G$ quasi-isometric to $H$ must
fit into a short exact sequence:
$$0\longrightarrow F\longrightarrow G\longrightarrow \Gamma
\longrightarrow 0$$ where $F$ is a finite group and $\Gamma\subset
Isom(\tilde X)$ (see Lemma 3.1).  The theorem follows.
\end{Prf}

For other recent results on the quasi-isometry behavior of
graphs of groups, we refer to Papasoglu [14], Papasoglu-Whyte [15],
and Mosher-Sageev-Whyte [13].

\section{Concluding remarks.}

We note that the only place in our arguments where we use the assumption
that $n=3$ is in the proof of the strong Jordan separation theorem.  
More specifically,
we make use of the fact that the Schoenflies theorem holds in dimension $2$.
Of course, this approach 
fails in higher dimension, as there are examples of wild
embeddings of spheres in all dimensions $\geq 3$.  Nevertheless, we still
believe that the conjecture put forth in the introduction holds true
(and in fact, that the strong Jordan separation theorem also holds in 
higher dimension).  

A more general question would be to determine which hyperbolic P-manifolds
exhibit rigidity.  One would need some sort of hypotheses, as the following
example shows:

\vskip 10pt

\noindent {\bf Example:}  Let $X_i$ ($1\leq i\leq 3$) be simple, thick, 
hyperbolic P-manifolds of dimension $3$.  In each $X_i$, let $Y_i\subset X_i$
be one of the surfaces in the $2$-dimensional strata, and let $\gamma _i$
be a simple closed geodesic in each $Y_i$.  We now propose to build a 
thick, hyperbolic P-manifold which is {\it not} rigid.  Let $G$ be a 
complete graph on four vertices, and $T_i$ be three triangles in $G$.
Assign a length to each edge in such a way that the triangles $T_i$ have
length equal to the corresponding $\gamma_i$.

Note that given isometries from $\gamma_i$ to $T_i$, we can form a
thick, hyperbolic P-manifold by gluing the $X_i$ to $G$ by identifying the
$\gamma_i$ with the $T_i$.  Furthermore, as long as the gluing isometries
are homotopic, the resulting P-manifolds will all have isomorphic fundamental
group.  So provided we can find two different gluings which yield 
non-isometric P-manifolds, we will have exhibited non-rigidity.  

To do this, fix the gluing of $X_2$, $X_3$, and `rotate' the gluing map for 
the $X_1$.  Note that, in the surface $Y_1$, there are countably many geodesic 
segments which are perpendicular to $\gamma_1$ and intersect $\gamma _1$ 
precisely at their endpoints.  
By rotating the gluing map suitably,
we can ensure that one of the resulting P-manifolds has such a geodesic
segment in $Y_1$ emanating from a vertex of the triangle $T_1$, whereas 
another one of the resulting P-manifolds does not.  It is now clear that
these two P-manifolds cannot be isometric, despite the fact that they
have isomorphic fundamental groups.  Observe
that these examples will have a non-trivial one-dimensional strata (namely
the graph $G$), so do not satisfy the simplicity hypothesis of this paper.

\vskip 10pt

Note that if there is no 
1-dimensional strata, one cannot use the `rotation trick' to get 
counterexamples (since the isometry group of compact hyperbolic manifolds
is finite if the dimension is at least two).  Perhaps the following is
reasonable:

\vskip 5pt

\noindent {\bf Question:} Is every hyperbolic P-manifold with empty 1-dimensional
strata Mostow rigid?

\vskip 5pt

Other interesting questions arise from trying to further understand the (full)
group of isometries $Isom(\tilde X)$ of the universal cover of a hyperbolic 
P-manifold $X$.  
These groups will be discrete, and exhibit behavior which one would 
expect to be between that of tree lattices (particularly if all the chambers
are isometric), and that of lattices in $SO(n,1)$.

Finally, we point out that we can define {\it negatively curved} P-manifolds
by allowing metrics of negative curvature (with totally geodesic boundary) 
on each chamber.  In this setting, the proof of Proposition 2.2 still holds (with 
appropriate modifications in the argument for Claim 2), and while there is
no hope of a Mostow type rigidity, one can still consider the quasi-isometry
question for the resulting groups.  The main stumbling point lies in the 
(3-dimensional case of the) following:

\vskip 5pt

\noindent {\bf Question:} Let $M^n$ be a compact, negatively curved, $n$-dimensional
manifold with non-empty, totally geodesic boundary, and let $\tilde M^n$ be
it's universal cover with the induced metric.  Is every quasi-isometry of 
$\tilde M^n$ a bounded distance from an isometry?

\vskip 5pt

\section{Bibliography.}

\vskip 10pt

\noindent [1] Bing, R.H.  {\it The geometric topology of 3-manifolds}.
American Math. Soc., Providence, RI, 1983.

\vskip 5pt

\noindent [2] Borel, A. \& Harish-Chandra.
{\it Arithmetic subgroups of algebraic groups}.
Ann. of Math. (2) 75 (1962), pp. 485--535.

\vskip 5pt

\noindent [3] Bourdon, M.
{\it Structure conforme au bord et flot geodesiques d'un $CAT(-1)$-espace}.
Enseign. Math. (2) 41 (1995), pp. 63--102.

\vskip 5pt

\noindent [4] Bourdon, M. \& Pajot, H. {\it Rigidity of quasi-isometries for
some hyperbolic buildings}.  Comment. Math. Helv. 75 (2000), pp. 701--736.

\vskip 5pt

\noindent [5] Bridson, M.R. \& Haefliger, A. {\it Metric spaces of 
non-positive curvature}. Springer-Verlag, Berlin, 1999.

\vskip 5pt

\noindent [6] Farb, B.  {\it The quasi-isometry classification of 
lattices in semisimple Lie groups}.  Math. Res. Lett. 4 (1997),
pp. 705-717.

\vskip 5pt

\noindent [7] Frigerio, R.  {\it Hyperbolic manifolds with geodesic
boundary which are determined by their fundamental group}.
\noindent preprint, http://front.math.ucdavis.edu/math.GT/0306398.

\vskip 5pt

\noindent [8] Forester, M.
{\it Deformation and rigidity of simplicial group actions on trees}.
Geom. Topol. 6 (2002), pp. 219--267.

\vskip 5pt

\noindent [9] Gromov, M. \& Piatetski-Shapiro, I.
{\it Nonarithmetic groups in {L}obachevsky spaces}.
 Inst. Hautes \'Etudes Sci. Publ. Math. 66 (1988), pp. 93--103.

\vskip 5pt

\noindent [10] Gromov, M. \& Thurston, W.  {\it Pinching constants for 
hyperbolic manifolds}.  Invent. Math. 89 (1987), pp. 1--12.

\vskip 5pt

\noindent [11] Guirardel, V.  {\it A very short proof of Forester's
rigidity result}.  Geom. Topol. 7 (2003), pp. 321-328.

\vskip 5pt

\noindent [12] Kapovich, M. \& Kleiner, B.  {\it Hyperbolic groups with
low-dimensional boundary}.  Ann. Sci. \'Ecole Norm. Sup. (4) 33 (2000),
pp. 647--669.

\vskip 5pt

\noindent [13] Mosher, L., Sageev, M. \& Whyte, K.  {\it Quasi-actions
on trees: research announcement}.
\noindent preprint, http://front.math.ucdavis.edu/math.GR/0005210.

\vskip 5pt

\noindent [14] Papasoglu, P.
{\it Group splittings and asymptotic topology}.
 
\noindent preprint, http://front.math.ucdavis.edu/math.GR/0201312.

\vskip 5pt

\noindent [15] Papasoglu, P. \& Whyte, K.  {\it Quasi-isometries between
groups with infinitely many ends}.  Comment. Math. Helv. 77 (2002),
pp. 133--144.

\vskip 5pt

\noindent [16] Paulin, F.  {\it Actions de groupes sur les arbres}, in
Seminaire Bourbaki, Vol. 1995/96.  Asterisque 241 (1997), pp. 97--137.

\vskip 5pt

\noindent [17] Sela, Z.  {\it Structure and rigidity in (Gromov) hyperbolic
groups and discrete groups in rank $1$ Lie groups. II.}  Geom. Funct. Anal.
7 (1997), pp. 561--593.

\vskip 5pt

\noindent [18] Sela, Z.  {\it Endomorphisms of hyperbolic groups. I. The Hopf
property}.  Topology 38 (1999), pp. 301--321.

\vskip 5pt

\noindent [19] Serre, J.P.  {\it Trees}.  Springer-Verlag, Berlin, 1980.

\vskip 5pt

\noindent [20] Sierpinski, W.  {\it Un th\'eor\`eme sur les continus}.  Tohoku
Math. Journ. 13 (1918), pp. 300--303.

\vskip 5pt

\noindent [21] Vinberg, \`E.B.  {\it Discrete reflection groups in Lobachevsky
spaces}, in Proceedings of the I.C.M. (Warsaw 1983), pp. 593--601.

\end{document}